\documentclass[letterpaper, 10 pt, conference]{ieeeconf}  

\pdfoutput=1

\IEEEoverridecommandlockouts                              

\usepackage[pdftex]{graphicx}
\usepackage{cite}
\usepackage{tikz}
\usepackage{nomencl}
\makenomenclature

\nomenclature{$f_i(t)$}{The friction force of $i$-th block for $i=1$ or $2$ at time $t$}%
\nomenclature{$f_u$}{The upper bound of the control force $f$: $\lvert f\rvert\leq f_u$}%
\nomenclature{$x_i(t)$}{The absolute position of $i$-th block for $i=1$ or $2$ at time $t$}%
\nomenclature{$u(t)$}{The center of mass velocity at time $t$: $u(t)=\frac{\dot{x_2}(t)+\dot{x_1}(t)}{2}$}%
\nomenclature{$v(t)$}{the difference of the velocities of two parts: $v(t)=\frac{\dot{x_2}(t)-\dot{x_1}(t)}{2}$ at time $t$}
\nomenclature{$d(t)$}{The position difference of two parts: $d(t)=x_2(t)-x_1(t)$ at time $t$}%
\nomenclature{$f_{bw}$}{The friction force in the backward motion}%
\nomenclature{$f_{fw}$}{The friction force in the forward motion}%
\nomenclature{$\rho$}{The ratio of two friction coefficient: $\rho=\frac{f_{fw}}{f_{bw}}$}%
\nomenclature{$\beta$}{The dimensionless constant $\beta=\frac{f_{bw}-f_{fw}}{2}$}%
\nomenclature{$\alpha$}{The dimensionless constant $\alpha=\frac{f_{fw}+f_{bw}}{2}$}%
\nomenclature{$m_i$}{The mass of $i$-th block for $i=1$ or $2$}%
\nomenclature{$T$}{The period of the motion}%
\nomenclature{$E(T)$}{The maximum size change in one period, \\$E(T)=\max_{t\in[t_1,t_1+T]}(d(t))-\min_{t\in[t_1,t_1+T]}(d(t))$}%
\nomenclature{$\eta$}{The ratio of maximum and minimum actuator force $\frac{f_u}{f_{bw}}$}%

\usetikzlibrary{calc,patterns,decorations.pathmorphing,decorations.markings}
\usepackage[cmex10]{amsmath}
\usepackage{amssymb}

  \usepackage[caption=false,font=footnotesize]{subfig}

\title{\LARGE \bf
Optimal periodic locomotion for a two piece worm with an asymmetric dry friction model
}

\author{Nak-seung Patrick Hyun$^{1}$ and Erik Verriest$^{2}$
\thanks{*This work was supported by a seed grant(RIM 3661304) from the Center for Robotics and Intelligent Machines at Georgia Institute of Technology}
\thanks{$^{1}$School of Electrical and Computer Engineering Georgia Institute of Technology
Atlanta, Georgia 30332--0250
        {\tt\small nhyun3@gatech.edu}}%
\thanks{$^{2}$School of Electrical and Computer Engineering Georgia Institute of Technology
Atlanta, Georgia 30332--0250
        {\tt\small erik.verriest@ece.gatech.edu}}%
}

\begin{document}

\maketitle
\thispagestyle{empty}
\pagestyle{empty}

\begin{abstract}
This paper solves the optimization problem for a simplified one-dimensional worm model when the friction force depends on the direction of the motion. The motion of the worm is controlled by the actuator force $f(t)$ which is assumed to be piecewise continuous and always generates the same force in the opposite directions. The paper derives the necessary condition for the force which maximizes the average velocity or minimizes the power over a unit distance. The maximum excursion of the worm body and the force are bounded. A simulation is given at the end of the paper.
\end{abstract}

%
\IEEEpeerreviewmaketitle

\section{Introduction}
The locomotion for legless animals, such as a snake and a worm, has many interesting features. One of the most appealing characteristics of these types of animal is that the periodic change in its shape generates the motion in a certain direction.  Scales on the animal's skin create different friction forces depending on the direction of the movement. Due to this asymmetric friction, these types of animal, such as snake or worm, can move forward. As a result, many researchers have focused on analyzing the periodic kinematics and dynamics of the legless animal locomotion in an asymmetric friction model, Verriest \cite{verriest2008locomotion} and Zimmermann \cite{zimmermann2003approach}. 

Even for a simplified one-dimensional worm model it is complicated to analyze its full dynamics. In the past years, a lot of research analyzed the one-dimensional toy problem in depth, Zimmermann\cite{Zimmermann2009book} and Chernousko\cite{chernous2002optimum}. Especially, the book \cite{Zimmermann2009book} has successfully characterized the kinematics and dynamics of the one-dimensional worm problem. 

Several researchers, Figurina\cite{figurina2007optimum} and Bolotnik\cite{bolotnik2007optimum}, have attempted to solve the optimization problem for a similar one-dimensional worm problem. Chernousko\cite{chernous2002optimum} has solved the optimal problem when the actuator force is modeled to be piecewise constant and for the case when the infinite force is available. He found the optimal control which maximizes the average velocity. Chernousko\cite{chernous2011optimum} also optimizes the energy consumed over a unit distance for all piecewise constant forces. Here, we attempt to optimize a similar performance index for an extended admissible control which includes all the piecewise continuous forces. In addition, we pose the realistic constraint that the force and the size of the worm are bounded.

In this paper, we first list the variables that are commonly used in this paper, and define the model which includes the friction model and the actuator model. Second, we suggest two performance indices to be optimized and state the problem formally. Third, we solve the problem analytically and find the necessary conditions for the optimal solution. Lastly, we provide one simulation example to understand the solution better.

\section{Modeling}

\begin{figure}[htp]
\centering
\begin{tikzpicture}[every node/.style={draw,outer sep=0pt,thick}]
	\tikzstyle{spring}=[thick,decorate,decoration={zigzag,pre length=0.3cm,post length=0.3cm,segment length=6}]
	\tikzstyle{dampener}=[thick,decorate,decoration={markings,mark connection node=dmp,mark=at position 0.5 with 
	  {
	    \node (dmp) [thick,inner sep=0pt,transform shape,rotate=-90,minimum width=15pt,minimum height=3pt,draw=none] {};
	    \draw [thick] ($(dmp.north east)+(2pt,0)$) -- (dmp.south east) -- (dmp.south west) -- ($(dmp.north west)+(2pt,0)$);
	    \draw [thick] ($(dmp.north)+(0,-5pt)$) -- ($(dmp.north)+(0,5pt)$);
	  }
	}]
	\tikzstyle{ground}=[fill,pattern=north east lines,draw=none,minimum width=0.75cm,minimum height=0.3cm]
	
	\node (M2) [minimum width=1.5cm, minimum height=1cm] {$m_2$};
	\node (M3) at (M2.west) [xshift=-3cm,minimum width=1.5cm, minimum height=1cm] {$m_1$};
	\draw [dampener] (M3.-0) to ($(M2.west)!(M3.-0)!(M2.west)$) ;
	\draw [->,thick] (M3.north west) ++(3.45,0.6cm) -- node[draw=none,auto,swap]{$Motion$} +(2cm,0);
	\draw [<-,thick] (M3.north east) ++(-1.80,0.6cm) -- node[draw=none,auto,swap]{$Motion$} +(2cm,0);
	\draw [<->,thick] (M3.south east) ++(0.55,0.1cm) -- node[draw=none,auto,swap]{$f$} +(1.2cm,0);
	\draw [<-,thick] (M3.south east) ++(2.65,-0.2cm) -- node[draw=none,auto,swap]{$f_{FW}$} +(0.5cm,0);
	\draw [->,thick] (M3.south east) ++(-0.8,-0.2cm) -- node[draw=none,auto,swap]{$f_{BW}$} +(0.5cm,0);
	\node (wall1) [ground, minimum width=7cm,xshift=-1.8cm,yshift=-1.5cm]{};
	\draw [->,thick] (wall1.south west) ++(-0.00,-0.20cm) -- node[draw=none,auto,swap]{$x(t)$} +(7cm,0);
\end{tikzpicture}
\caption{One dimensional two piece worm model.}
\label{fig_model}
\end{figure}
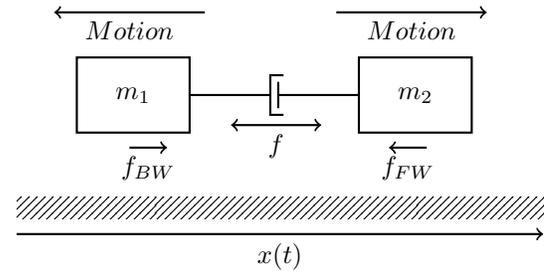
\printnomenclature	
\subsection{Configuration space}
Consider a simplified one-dimensional worm consisting of two parts, head and tail. These parts are connected by the actuator which exerts the same force $f(t)$ at time $t$ to both parts (Fig.\ref{fig_model}). This actuator introduces energy to system, which converts into kinetic energy for each parts and compensates the work done by the friction force. Since the worm contains two independent components, the dynamics of motion can be written as two uncoupled second order nonlinear differential equation. The nonlinearity comes from the differential friction which depends on the direction of the movement. More detail on the friction will be shown in the next subsection. 
\begin{equation}
\begin{array}{rl}
-f(t)+f_1(t)=m_1\ddot{x_1}(t)\\
f(t)+f_2(t) =m_2\ddot{x_2}(t)
\end{array} 
\label{eqn_dyn}
\end{equation}
In this paper, we define a shape space, $\mathbb{R}^+$, as an extend of the worm, $d(t)$. The shape can be changed by switching the sign of the force $f$ as the worm stretches or contracts its body. With this shape space, the original space can be prjected to the three dimensional space, where its trajectory $\Theta=\{(d(t),\dot{x_1}(t),\dot{x_2}(t))\subset \mathbb{R}^3| \forall t\in \mathbb{R}^+\}$ lies in $\mathbb{R}^3$, assuming that dynamics starts at $t=0$ and appropriate initial conditions are given. Using the change of basis method, we now define the configuration space in $\mathbb{R}^3$ with its homeomorphic trajectory $\Omega=\{(d(t),v(t),u(t))\subset \mathbb{R}^3|$ for $\forall t\in \mathbb{R}^+\}$. This paper analyzes the trajectory, $\Omega$ of the worm.
\subsection{Friction model}
Our friction model is an asymmetric dry friction which depends on the sign of the velocity but is independent of the relative speed. The asymmetry in the friction will result in moving the worm in a certain direction. More concretely, let $f$ be an asymmetric Coulomb friction which is an approximate model for dry friction. In addition, let $m=m_1=m_2$ so that each part experiences the same friction.
\begin{equation}
\begin{array}{rl}
f_i(t) = \left\{ \begin{array}{rl}
 -f_{fw} &\mbox{ if $\dot{x_i}(t)>0$} \\
  f_{bw} &\mbox{ if $\dot{x_i}(t)<0$} \\
  f_0		 &\mbox{ if $\dot{x_i}(t)=0$}
\end{array} \right.
\end{array}
\label{eqn_friction}
\end{equation}
where $f_{fw},f_{bw},f_0>0$ are constant and $i=1$ or $2$. The friction is modeled to satisfy $f_{bw}\geq f_{fw}$ assuming that the worm moves in forward direction. As a result of asymmetry, the center of mass will move forward. In other words, the third component in $\Omega$: $u(t)$, will always be nonnegative for any $t\in\mathbb{R}^+$. 

\subsection{Periodic motion in configuration space}
\begin{figure}[!ht]
\centering
\subfloat[$(v,u)$ plane]{\includegraphics[width=2.9in]{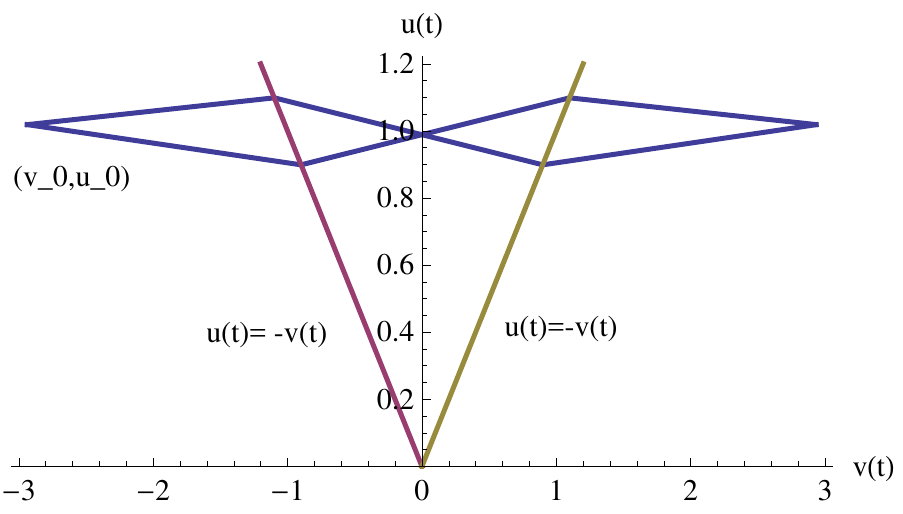}
\label{fig_first_case}}
\vfil
\subfloat[$(d,v)$ plane]{\includegraphics[width=2.9in]{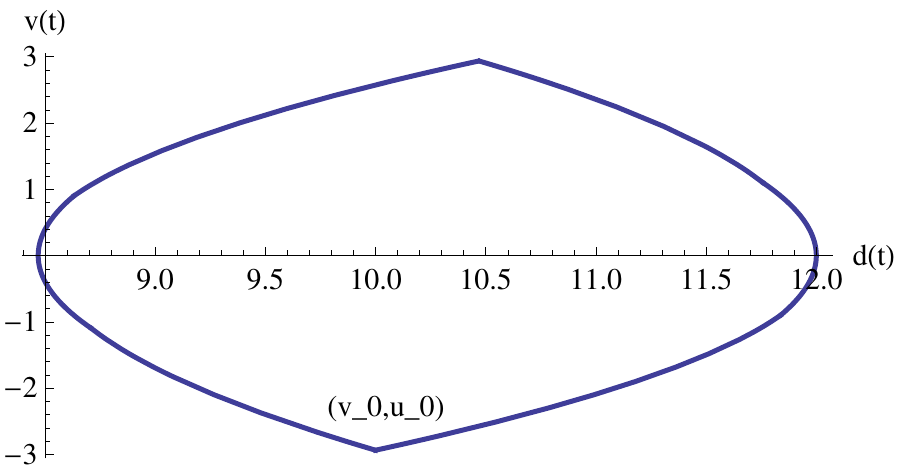}
\label{fig_second_case}}
\caption{Simulation results for $\rho=0.5, d_0=10$ and $F=5$}
\label{fig_sim}
\end{figure}

There are eight cases of dynamics depending on the sign of the velocities and the sign of the force acting on each part. 
\begin{table}[!ht]
\renewcommand{\arraystretch}{1.3}
\caption{Eight modes of dynamics.}
\label{table_example}
\centering
\begin{tabular}{c||c||c||c||c||c||c||c||c}
\hline
\bfseries CASE &\bfseries 1 & \bfseries 2 & \bfseries 3 & \bfseries 4 & \bfseries 5 & \bfseries 6 & \bfseries 7 & \bfseries 8 \\
\hline\hline
$f(t)$ & $+$ & $+$ & $+$ & $-$ & $-$ & $-$& $+$ & $-$ \\
\hline
$\dot{x_1}(t)$ & $+$ & $+$ & $-$ & $-$ & $+$ & $+$& $-$ & $-$\\
\hline
$\dot{x_2}(t)$ & $-$ & $+$ & $+$ & $+$ & $+$ & $-$& $-$ & $-$\\
\hline
$Duration$ & $T_1$ & $T_2$ & $T_3$ & $T_4$ & $T_5$ & $T_6$& $T_7$ & $T_8$\\
\hline
\end{tabular}
\label{tbl_cases}
\end{table}
TABLE (\ref{tbl_cases}) shows the resulting eight possible cases of dynamics. The table also contains the duration of each case indicated by $T_i$ for each case $i$. Observe that CASE 7 and CASE 8 cannot occur since $u(t)$ is always nonnegative. It suffices to analyze the motion for the six cases. For conveninece, $m_i=1$ for all $i$. Using the friction model, Eqn(\ref{eqn_friction}), and the dynamics, Eqn(\ref{eqn_dyn}), we derive Eqn(\ref{eqn_cases_v}) and Eqn(\ref{eqn_cases_u}). With initial conditions $(v(t_i),u(t_i))=(v_i,u_i)$ for each $i$ case, the Eqn$(\ref{eqn_dyn})$ integrates to,

\begin{IEEEeqnarray}{c}v(t)=\left\{ \begin{array}{rl}
v_1+\int_{t_1}^{t_1+T_1}{|f(t)|dt}+\alpha(t-t_1) &\mbox{CASE 1}\\
v_2+\int_{t_2}^{t_2+T_2}{|f(t)|dt} &\mbox{CASE 2}\\
v_3+\int_{t_3}^{t_3+T_3}{|f(t)|dt}-\alpha(t-t_3) &\mbox{CASE 3}\\
v_4-\int_{t_4}^{t_4+T_4}{|f(t)|dt}-\alpha(t-t_4) &\mbox{CASE 4}\\
v_5-\int_{t_5}^{t_5+T_5}{|f(t)|dt} &\mbox{CASE 5}\\
v_6-\int_{t_6}^{t_6+T_6}{|f(t)|dt}+\alpha(t-t_6) &\mbox{CASE 6}\\
\end{array} \right.
\label{eqn_cases_v}
\\
u(t)=\left\{ \begin{array}{rl}
u_1+\beta(t-t_1) &\mbox{CASE 1}\\
u_2-f_{fw}(t-t_2) &\mbox{CASE 2}\\
u_3+\beta(t-t_3) &\mbox{CASE 3}\\
u_4+\beta(t-t_4) &\mbox{CASE 4}\\
u_5-f_{fw}(t-t_5) &\mbox{CASE 5}\\
u_6+\beta(t-t_6) &\mbox{CASE 6}\\
\end{array} \right.
\label{eqn_cases_u}
\end{IEEEeqnarray}
Similarly, the dynamics of the configuration variable, $d(t)$, is solved in terms of $v(t)$. With initial conditions $d(t_i)=d_i$ and $v(t_i)=v_i$ as before for each $i$ case, $d(t)$ can be derived by integrating the Eqn$(\ref{eqn_cases_v})$ 
\begin{equation}
d(t)=2\int_{t_i}^{t_i+T_i}v(t)dt + d_i \mbox{ for any $t\in[t_i,t_i+T_i]$ }
\label{eqn_cases_d}
\end{equation}
for all $i$. In order to find a periodic control which generates a periodic motion in configuration space, it suffices to check if there exists a periodic $f$ defined as Eqn(\ref{eqn_dyn}) such that $(d(t),v(t),u(t))=(d(t+T),v(t+T),u(t+T))$. If this is true then the trajectory of $(d(t),v(t),u(t))$ follows an orbit in $\Omega$ space and so does its projection to a subspace of $\Omega$. 
For example, a periodic orbit of the $\Omega$ trajectory is shown in Fig.\ref{fig_sim}, which $(a)$ shows the projection to the $(v,u)$ plane and $(b)$ shows the projection to the $(v,d)$ plane. Here we chose $|f(t)|=5$ constant, $(v(t_1),u(t_1),(d_1))=(-3,1,10)$ and $\rho=0.5$. 

\subsection{Periodic actuator force}
Most research on the periodic motion of the one dimensional worm model designed an actuator which generates a simple harmonic force \cite{Zimmermann2009book},\cite{fang2012dynamics} and analyzed its periodic dynamics. However, in this paper, the worm is assumed to generate an arbitrary periodic piecewise continuous force. This actuator force is modeled as follows, \begin{equation}
\begin{array}{rl}
&|f(t)|= F(t) \mbox{ $t\in[t_1,t_1+T]$ }
\\
&f(t+T)=f(t) \mbox{ for all $t\in \mathbb{R}^+$}
\end{array}
\label{eqn_actuator}
\end{equation}
where $T$ is the period and $F(t)$ is a positive piecewise continuous function. If $F(t)\geq f_{bw}$, the actuator can change the direction of motion of each part. Previous research \cite{miller1980existence},\cite{nistri1983periodic} shows a sufficient condition for the existence of periodic motion in certain nonlinear systems. Since the friction force and periodic actuator control are bounded, the nonlinearities in Eqn(\ref{eqn_dyn}) are bounded by $max_{t\in[t_1,t_1+T]}\{F(t)\}+f_{bw}$. By the  construction, both forces are continuous almost everywhere and also we have a freedom in choosing the sign of the control force. This satisfies the sufficient condition derived in \cite{nistri1983periodic}, and so there exists a periodic solution for Eqn(\ref{eqn_dyn}).

Here, we list four assumption for the actuator force design. First, we assume that the worm switches the sign of the actuator force once within a period. It is required to switch the sign of the actuator force at least once to pull back to the original state in the configuration space. This assumption shows that it tends to only extend its body as much as possible when $f(t)>0$, and contract its body back as the sign changes to negative, $f(t)<0$. In the presence of inertia, we assume that the switch happens when two parts move towards or move in the opposite direction in order to keep the maximum worm size. The physical constraint of the maximum allowable extension will appear in the next section. 

Second, following the design assumption, we enumerate a sequence of cases in TABLE (\ref{tbl_cases}) which the worm experiences. Suppose the force is initially set to positive at time $t_1$ as defined before, then $Case 1\to Case 2\to Case 3$ will be the motion sequence for $f(t)>0$ in $t\in[t_1,t_1+T_1+T_2+T_3]$. After the sign switches at $t_1+T_1+T_2+T_3$, the sequence changes to $Case 4\to Case 5\to Case 6$ for $f(t)<0$ in $t\in[t_1+T_1+T_2+T_3,t_1+T]$. This shows that $t_i=t_1+\sum_{k=1}^{i-1}T_k$ for all $i={2,\cdots,6}$.

Third, there exist a symmetry on the actuator force design between the positive and negative region. In other words, the initial state in $(v,u)$-plane for negative force is equal to $(v_4,u_4)=(-v_1,u_1)$. By using the definition of $(v, u)$, this shows that $(\dot{x_1}(t_4),\dot{x_2}(t_4))=(\dot{x_2}(t_1),\dot{x_1}(t_1))$. Therefore, by designing the negative force $f(t)=-f(t-\frac{T}{2})$, the dynamics in Eqn(\ref{eqn_dyn}) are the same except for the negative sign on each equation. This guarantees that such a force design will generate the periodic motion in $(v,d)$-space since the integral term in  Eqn(\ref{eqn_cases_v}) and (\ref{eqn_cases_d}) will cancel out and so 
\begin{IEEEeqnarray*}{c}
\begin{array}{rl}
v(t_1+T)&=v(t_1)
\\
d(t_1+T)&=d(t_1).
\end{array}
\end{IEEEeqnarray*}

Lastly, we define a relation between $T_1, T_2$ and $T_3$ in order to have $u(t_1)=u(t_1+T)$ which finalize the full periodic motion in the configuration space. From the above third assumption, we know that $T_1=T_4, T_2=T_5, T_3=T_6$ and so $\frac{T}{2}=T_1+T_2+T_3=T_4+T_5+T_6$. Here, we claim that if $T_2=\frac{1-\rho}{2\rho}(T_1+T_3)$ then $u(t_1+T)=u(t_1)$. By following the sequence defined above and Eqn(\ref{eqn_cases_u}), we have 
\begin{IEEEeqnarray*}{c}
\begin{array}{rl}
u(t_1+T)&=u_1+\beta(T_1+T_3+T_4+T_6)-f_{fw}(T_2+T_5)
\\
&=u_1+2\beta(T_1+T_3)-2f_f(T_2)
\\
&=u_1.
\end{array}
\end{IEEEeqnarray*}

In addition, we now have a useful set of equations. 
\begin{IEEEeqnarray}{c}
\begin{array}{rl}
&T_2=\frac{T}{2}(\frac{1+\rho}{1-\rho})
\label{eqn_T2}
\end{array}
\\
\begin{array}{rl}
&T_3=\frac{T}{2}(\frac{2\rho}{1+\rho})-T_1
\label{eqn_T3}
\end{array}
\end{IEEEeqnarray}
Hence, for given $T$ and $T_1$, any piecewise continuous function $F(t)$ for $t\in[t_1,t_1+\frac{T}{2}]$ can generate a periodic motion with period $T$ in the configuration space. 

\subsection{Physical constraints}

There are two physical constraints for the worm in this paper. One is the size of the worm and the other is the upper limit of the actuator force. The first one states that the difference between maximum and minimum excursion of the worm's body is limited to $E(T)\leq L$ and the actuator was designed to reach this bound, $L$. The other constraint is that there exists an upper limit, $f_u>0$, of the actuator force $F$ because of the physical limit of the worm's muscle. With these two constraints, we can derive $F$ which optimizes given performances such as the power over unit distance traveled and the speed of the center of the mass. Here, we define a set of admissible control forces by
\begin{equation}
\begin{array}{rl}
S_{T}=&\{F \in PC[t_1,t_1+\frac{T}{2}]:F(t)\in[f_{bw},f_u]
\\
&\mbox{for all $t\in[t_1,t_1+\frac{T}{2}]$}\}
\end{array}
\end{equation}
where $PC$ is a set of piecewise continuous function.

\section{Optimization}
\subsection{Performance index}
Research on the optimization of this one dimensional worm model was first suggested in \cite{chernous2002optimum}, and extended on \cite{chernous2011optimum}. The latter paper suggests two performance indices which we adopt in this paper. One is the power over the unit distance and the other is the average velocity over one period. Other papers in the optimization for one-dimensional worm problem, \cite{figurina2007optimum}\cite{bolotnik2007optimum}, also focused on maximizing the average velocity over one period.

Let Eqn(\ref{eqn_X_P}) be the total distance that the center of mass traveled and the average power exerted by the actuator, respectively. 
\begin{equation}
\begin{array}{rl}
&X=\int_{t_1}^{t_1+T}u(t)dt
\\
&P=\frac{1}{T}\int_{t_1}^{t_1+T}F(t)v(t)dt
\label{eqn_X_P}
\end{array}
\end{equation}
The two performance indeces in this paper are as follows, 
\begin{IEEEeqnarray}{c}
\begin{array}{rl}
V=\frac{X}{T}
\label{eqn_V_T}
\end{array}
\\
\begin{array}{rl}
P_u=\frac{P}{X}.
\label{eqn_P_U}
\end{array}
\end{IEEEeqnarray}
These performance indices are the average velocity over one period, $V$, and the average power over unit distance, $P_{u}$. 

\subsection{Optimal control problem}
Past research of the optimal control problem on one dimensional two body system have modeled the actuator force to be constant in each dynamic motion, \cite{chernous2002optimum} and \cite{chernous2011optimum}. In this paper, we not only include the constant model but also generalize to any piecewise continuous force model in each dynamics.

There are two optimal control problems in this paper for the fixed period $T$. The problem is stated for each performance index,
\begin{enumerate}
\item Find an optimal $F\in S_{T}$ for a given $T_1$ with the maximum excursion constraint, $E(T)=L$.
\item Find an optimal $F\in S_{T}$ with the maximum excursion limit, $E(T)=L$.
\end{enumerate}

\subsection{Solution to the first problem}
Pick $T_1$ such that $T_1\in [0,\frac{T}{2}-T_2]$, then do the following:
\begin{enumerate}
\item STEP 1: Compute the boundary condition for each dynamics and $E(T)$
\item STEP 2: Compute two performance index, $P_u$ and $V$ in Eqn(\ref{eqn_P_U}) and Eqn(\ref{eqn_V_T}).
\item STEP 3: Solve the optimization problem for a given $T_{min}$.
\item STEP 4: Find the necessary and sufficient conditions for the existence of solution.
\item STEP 5: Find the optimal $T_{min}$ which minimize the performance.
\end{enumerate}
Here are the details:
\subsubsection{STEP 1}
First, we compute the boundary condition for each dynamics. To do this, we need to define three additional functions, $H:[t_1,t_1+T_1]\to \mathbb{R}, G:[t_2,t_2+T_2]\to \mathbb{R}$ and $I:[t_3,t_3+T_3]\to \mathbb{R}$, which satisfy
\begin{IEEEeqnarray}{c}
H(t)=\int_{t_1}^{t}F(t)dt
\label{eqn_H}
\\
G(t)=\int_{t_2}^{t}F(t)dt
\label{eqn_G}
\\
I(t)=\int_{t_3}^{t}F(t)dt.
\label{eqn_I}
\end{IEEEeqnarray}
The sequence of motion and the TABLE {\ref{tbl_cases}} shows that $\dot{x_2}(t_2)=0$ and $\dot{x_1}(t_3)=0$. By the dynamics of ${x_i}(t)$ for $i\in\{1,2\}$, we derive the boundary conditions for $H$ and $G$. In addition, by using the third assumption of the actuator force, $v_4=-v_1$, we derive the boundary condition for $I$ as well. Here, we use  Eqn({\ref{eqn_cases_v}}), Eqn({\ref{eqn_T2}}) and Eqn(\ref{eqn_T3}) for computation.
\begin{IEEEeqnarray*}{c}
\begin{array}{rl}
\dot{x_2}(t_2)&=\int_{t_1}^{t_1+T_1}F(t)dt+f_{bw}T_1+(u_1+v_1)
\\
&=H(t_1+T_1)+f_{bw}T_1+(u_1+v_1)
\\
&=0
\label{eqn_H_bound}
\end{array}
\\
\begin{array}{rl}
\dot{x_1}(t_3)&=\int_{t_2}^{t_2+T_2}-F(t)dt-f_{fw}T_2+\dot{x_1}(t_2)
\\
&=-G(t_2+T_2)-f_{fw}T_2+(2u_1+2\beta T_1)
\\
&=0
\label{eqn_G_bound}
\end{array}
\\
\begin{array}{rl}
\dot{v}(t_4)&=\int_{t_3}^{t_3+T_3}F(t)dt-\alpha T_3+\beta T_1-f_{fw}T_2+u_1
\\
&=I(t_3+T_3)-f_{bw}T_1-\frac{Tf_{fw}}{1+\rho}+u_1
\\
&=-v_1.
\label{eqn_I_bound}
\end{array}
\end{IEEEeqnarray*}
By using Eqn(\ref{eqn_T2}) and the above equations, the boundary conditions for $H, G$ and $I$ are
\begin{IEEEeqnarray}{c}
H(t_1+T_1)=-f_{bw}T_1+c_1
\label{eqn_H_bound_final}
\\
G(t_2+T_2)=2\beta T_1+c_2
\label{eqn_G_bound_final}
\\
I(t_3+T_3)=-f_{bw}T_1+c_3
\label{eqn_I_bound_final}
\end{IEEEeqnarray}
where $c_1=-(u_1+v_1), c_2=2u_1-\frac{Tf_{fw}}{2}\frac{1-\rho}{1+\rho}$ and $c_3=-(u_1+v_1)+\frac{Tf_{fw}}{1+\rho}$ are constant.

Second, we need to compute the difference between maximum and minimum excursion, $E(T)$. To do this, we need to find the time when the size of the worm, $d(t)$, are at minimum and maximum. Since $v(t)<0$ for $Case 1$, the Eqn(\ref{eqn_cases_d}) is decreasing in $[t_1,t_1+T_1]$. However, we know that $v(t_2+T_2)>0$ and so there exist some $T_{min}\in[0,T_2]$ such that $v(t_2+T_{min})=0$. Similarly, $v(t)>0$ for $Case 4$ and so $d(t)$ is increasing in $[t_4,t_4+T_4]$ but $v(t_5+T_5)<0$. Therefore, there also exist $T_{max}\in[0,T_5]=[0,T_2]$ such that $v(t_5+T_{max})=0$. The symmetry assumption for the actuator force shows that $T_{min}=T_{max}$. From this we can compute $E(T)$ as follows.
\begin{IEEEeqnarray*}{c}
\begin{array}{rl}
E(T)&=2\int_{t_2+T_{min}}^{t_5+T_{max}}v(t)dt
\\
&=2(\int_{t_2+T_{min}}^{t_4}v(t)dt+\int_{t_4}^{t_5+T_{max}}v(t)dt)
\\
&=2(\int_{t_2+T_{min}}^{t_4}v(t)dt-\int_{t_1}^{t_2+T_{min}}v(t)dt)
\\
&=2(-\int_{t_1}^{t_2+T_{min}}v(t)dt+\int_{t_2+T_{min}}^{t_4}v(t)dt).
\label{eqn_excursion_compute}
\end{array}
\end{IEEEeqnarray*}
By applying Eqn(\ref{eqn_cases_v}) to the above equation, we get 
\begin{IEEEeqnarray*}{c}
\begin{array}{rl}
E(T)&=2(-(\int_{t_1}^{t_2}H(t)dt+\int_{t_1}^{t_2}(\alpha t+ v_1)dt)
\\
&-(\int_{t_2}^{t_2+T_{min}}G(t)dt+\int_{t_2}^{t_2+T_{min}}(-\alpha T_1-u_1)dt)
\\
&+(\int_{t_2+T_{min}}^{t_3}G(t)dt+\int_{t_2+T_{min}}^{t_3}(-\alpha T_1-u_1)dt)
\\
&+(\int_{t_3}^{t_4}I(t)dt+\int_{t_3}^{t_4}(-\alpha t +f_{bw}T_1+\beta T_2+u_1)dt))
\\
&=2h(T_1,u_1,v_1,T_{min})
\\
&+2(\int_{t_3}^{t_4}I(t)dt-\int_{t_1}^{t_2}H(t)dt)
\\
&-2(\int_{t_2}^{t_2+T_{min}}G(t)dt-\int_{t_2+T_{min}}^{t_3}G(t)dt)
\label{eqn_excursion_compute_final}
\end{array}
\end{IEEEeqnarray*}
where $2h(T_1,u_1,v_1,T_{min})$ collects the constant terms of the equation. By the constraint of this optimization problem, we want $E(T)=L$ and so finally we get,
\begin{IEEEeqnarray}{c}
\begin{array}{rl}
\int_{t_3}^{t_4}I(t)dt-\int_{t_1}^{t_2}H(t)dt&=\frac{L}{2}-h(T_1,u_1,v_1,T_{min})
\\
&+\int_{t_2}^{t_2+T_{min}}G(t)dt
\\
&-\int_{t_2+T_{min}}^{t_3}G(t)dt
\label{eqn_excursion_compute_final_L}
\end{array}
\end{IEEEeqnarray}

In addition, $v(t_2+T_{min})=0$ shows that $\int_{t_2}^{t_2+T}F(t)dt-\beta T_1-u_1=0$ and so we have
\begin{equation}
G(t_2+T_{min})=\beta T_1+u_1.
\label{eqn_G_internal}
\end{equation}

\subsubsection{STEP 2}
Now, we compute the performance index. We start with the average velocity over one period defined in Eqn(\ref{eqn_V_T}). To do this, we compute $X$ first. By using the Eqn(\ref{eqn_cases_u}), Eqn(\ref{eqn_X_P}), Eqn(\ref{eqn_T2}) and Eqn(\ref{eqn_T3}), we get
\begin{IEEEeqnarray*}{c}
\begin{array}{rl}
X=\int_{t_1}^{t_1+T}u(t)dt&=2\int_{t_1}^{t_4}u(t)dt
\\
&=T(u_1+\beta T_1-\frac{Tf_{fw}\beta}{4\alpha}).
\end{array}
\end{IEEEeqnarray*}
By simply dividing by $T$, we get the average velocity performance index,
\begin{IEEEeqnarray}{c}
V=u_1+\beta T_1-\frac{Tf_{fw}\beta}{4\alpha}
\label{eqn_V_final}
\end{IEEEeqnarray}

Next, we compute the second performance index, $P_u$. Since $X$ was constant for fixed $T_1$, it is enough to optimize the power over a unit distance. The symmetry assumption in the actuator force design shows that the total work, $W_{T,T_1}$, is equal to two times of the work done in the half of a period, $2W_{T/2,T_1}$. Therefore, it suffices to compute the work done in the half of a period.
\begin{IEEEeqnarray}{c}
\begin{array}{rl}
W_{T/2,T_1}&=\int_{t_1}^{t_1+\frac{T}{2}}F(t)v(t)dt
\\
&=W_1+W_2+W_3
\label{eqn_work}
\end{array}
\end{IEEEeqnarray}
where $W_i=\int_{t_i}^{t_i+T_i}F(t)v(t)dt$ are the work done in the $i-$th $Case$ for $i\in\{1,2,3\}$. By using the fundamental theorem of calculus, we can show that $H^\prime=F$ for a.e $t\in(t_1,t_1+T_1)$, $G^\prime=F$ for a.e $t\in(t_2,t_2+T_2)$ and $I^\prime=F$ for a.e $t\in(t_3,t_3+T_3)$, where a.e stands for almost everywhere. This shows that 
\begin{IEEEeqnarray*}{c}
\begin{array}{rl}
\int_{t_1}^{t_1+T_1}F(t)H(t)dt&=\int_{H(0)}^{H(t_1+T_1)}HdH
\\
&=\frac{H(t_1+T_1)^2}{2}
\end{array}
\\
\begin{array}{rl}
\int_{t_2}^{t_2+T_2}F(t)G(t)dt&=\int_{G(0)}^{G(t_1+T_1)}GdG
\\
&=\frac{G(t_2+T_2)^2}{2}
\end{array}
\\
\begin{array}{rl}
\int_{t_3}^{t_3+T_3}F(t)I(t)dt&=\int_{I(0)}^{I(t_1+T_1)}IdI
\\
&=\frac{I(t_3+T_3)^2}{2}.
\end{array}
\end{IEEEeqnarray*}
By using the above equations and Eqn(\ref{eqn_cases_v}), we can compute each $W_i(T_i)$ for all $i$. Here, we use integration by parts to simplify the equations, 
\begin{IEEEeqnarray}{c}
\begin{array}{rl}
W_1&=\int_{t_1}^{t_1+T_1}F(t)v(t)dt
\\
&=\frac{H(t_1+T_1)^2}{2}+\int_{t_1}^{t_1+T_1}(F(t)\alpha t +v_1F(t))dt
\\
&=\frac{H(t_1+T_1)^2}{2}+H(t_1+T_1)(\alpha T_1+v_1)
\\
&-\alpha\int_{t_1}^{t_1+T_1}H(t)dt
\\
&=W_{1const}-\alpha\int_{t_1}^{t_1+T_1}H(t)dt
\end{array}
\\
\begin{array}{rl}
W_2&=\int_{t_2}^{t_2+T_2}F(t)v(t)dt
\\
&=\frac{G(t_2+T_2)^2}{2}-(\beta T_1+u_1)\int_{t_2}^{t_2+T_2}F(t)dt
\\
&=\frac{G(t_2+T_2)^2}{2}-(\beta T_1+u_1)G(t_2+T_2)
\\
&=W_{2const}
\end{array}
\\
\begin{array}{rl}
W_3&=\int_{t_3}^{t_3+T_3}F(t)v(t)dt
\\
&=\frac{I(t_3+T_3)^2}{2}-(I(t_3+T_3)+v_1)I(t_3+T_3)
\\
&+\alpha\int_{t_3}^{t_3+T_3}I(t)dt
\\
&=W_{3const}+\alpha\int_{t_3}^{t_3+T_3}I(t)dt.
\end{array}
\end{IEEEeqnarray}
Let $W_{const}(T_1,T,u_1,v_1)=W_{1const}+W_{2const}+W_{3const}$ and by substituting above equations to Eqn(\ref{eqn_work}), we get 
\begin{IEEEeqnarray}{c}
\begin{array}{rl}
W_{T,T_1}&=2W_{\frac{T}{2},T_1}
\\
&=2W_{const}(T_1,T,u_1,v_1)
\\
&+2\alpha(\int_{t_3}^{t_3+T_3}I(t)dt-\int_{t_1}^{t_1+T_1}H(t)dt).
\end{array}
\end{IEEEeqnarray} 
By substituting the constraint, Eqn(\ref{eqn_excursion_compute_final_L}), we conclude that 
\begin{IEEEeqnarray}{c}
\begin{array}{rl}
W_{T,T_1}&=2W_{const}(T_1,T,u_1,v_1)
\\
&+\alpha L-2\alpha h(T_1,u_1,v_1,T_{min})
\\
&+2\alpha\int_{t_2}^{t_2+T_{min}}G(t)dt
\\
&-2\alpha\int_{t_2+T_min}^{t_3}G(t)dt.
\end{array}
\label{eqn_work_final_L}
\end{IEEEeqnarray}

\subsubsection{STEP 3}
Since $V$ did not depend on the choice of $F\in S_{T,T_1}$ for fixed $T_1$, there is nothing to optimize. For the second performance index, we summarize the problem as follows. For given $T_1$ and $L$, we want to find $F\in S_{T,T_1}$ such that minimize the total work, $W_{T,T_1}$ in Eqn(\ref{eqn_work_final_L}), and satisfies the boundary conditions in Eqn(\ref{eqn_H_bound_final}), Eqn(\ref{eqn_G_bound_final}), Eqn(\ref{eqn_I_bound_final}) and the excursion constraint in Eqn(\ref{eqn_excursion_compute_final_L}). Observe that the only function that changes the total work is $G$ since other term remain constant for a given $T_1$ and $T_{min}$. Therefore, it is free to choose $H$ and $I$ within the admissible control that satisfies all the constraints. Hence, the optimization problem reduces to find the optimal $G$ which optimizes $\min_{G\in PC[t_2,t_2+T_2]}(\int_{t_2}^{t_2+T_{min}}G(t)dt-\int_{t_2+T_{min}}^{t_3}G(t)dt)$ with the given boundary constraint. 

We regard $G(t)$ as a composition with two different dynamics
\begin{equation*}
\begin{array}{rl}
G(t) = \left\{ \begin{array}{rl}
 F(t) &\mbox{ if $t\in[t_2,t_2+T_{min}]$} \\
  -F(t) &\mbox{ if $t\in(t_2+T_{min},t_2+T_2]$}.
\end{array} \right.
\end{array}
\label{eqn_G_ham}
\end{equation*}
Accordingly, we construct two Hamiltonians,
\begin{equation*}
\begin{array}{rl}
H^1(t)=G(t)+\lambda_1(t)F(t)
\\
H^2(t)=G(t)-\lambda_2(t)F(t).
\end{array}
\end{equation*}
Since $F(t)\in[f_{bw},f_u]$ bounded, by applying the Pontryagain minimum principle(\cite{BrysonHo69}), we find the optimal control
\begin{equation*}
\begin{array}{rl}
F_1(t) = \left\{ \begin{array}{rl}
 f_{bw} &\mbox{ if $\lambda_1(t)>0$} \\
 f_u &\mbox{ if $\lambda_1(t)<0$}
\end{array} \right.
\end{array}
\end{equation*}
for $H^{1}$ Hamiltonian. Computing the Euler-Lagrange equation on $H^1(t)$ gives the differential equation for $\lambda_1(t)$ equation.
\begin{equation*}
\dot{\lambda_1}(t) = -\frac{\partial H^1}{\partial G}=-1
\end{equation*} 
which is solved by $\lambda_1(t)=-t+\tau_1$ for some $\tau_1$ constant. By using the internal point constraint, Eqn(\ref{eqn_G_internal}), we get 
\begin{IEEEeqnarray*}{c}
\begin{array}{rl}
G(t_2+T_{min})&=\int_{t_2}^{t_2+T_{min}}F(t)dt
\\
&=f_{bw}\tau_1+(T_{min}-\tau_1)f_u
\\
&=\beta T_1+u_1
\end{array}
\end{IEEEeqnarray*}
and solving for $\tau_1$ gives
\begin{IEEEeqnarray*}{c}
\begin{array}{rl}
\tau_1&=\frac{T_{min}f_u-(\beta T_1+u_1)}{f_u-f_{bw}}
\\
&=\frac{f_u-\frac{(\beta T_1+u_1)}{T_{min}}}{f_u-f_{bw}}T_{min}.
\end{array}
\end{IEEEeqnarray*}
Since $\tau_1\in[0,T_{min}]$, if $T_{min}\leq\frac{(\beta T_1+u_1)}{f_{bw}}$, the optimal solution exist and the solution is
\begin{equation*}
\begin{array}{rl}
F_1(t) = \left\{ \begin{array}{rl}
 f_{bw} &\mbox{ if $t\in[t_2,t_2+\tau_1]$} \\
 f_u &\mbox{ if $t\in[t_2+\tau_1,t_2+T_{min}]$}.
\end{array} \right.
\end{array}
\end{equation*}
Similarly, we solve for the second Hamiltonian system. By applying the Pontryagin minimum principle, we find the optimal control
\begin{equation*}
\begin{array}{rl}
F_2(t) = \left\{ \begin{array}{rl}
 f_u &\mbox{ if $\lambda_2(t)>0$} \\
 f_{bw} &\mbox{ if $\lambda_2(t)<0$}
\end{array} \right.
\end{array}
\end{equation*}
and the Euler-Lagrange equation gives 
\begin{equation*}
\dot{\lambda_2}(t) = -\frac{\partial H^2}{\partial G}=-1
\end{equation*} 
which is  solved by $\lambda_2(t)=-t+\tau_2$ for some $\tau_2$ constant. Now, by using the final state constraint, Eqn(\ref{eqn_G_bound_final}), we compute
\begin{IEEEeqnarray*}{c}
\begin{array}{rl}
G(t_2+T_2)&=\int_{t_2}^{t_2+T_{min}}F(t)dt+G(t_2+T_{min})
\\
&=f_{u}\tau_2+(T_2-(T_{min}+\tau_2))f_{bw}+\beta T_1+u_1
\\
&=2\beta T_1+c_2
\end{array}
\end{IEEEeqnarray*}
and solving for $\tau_2$ gives 
\begin{IEEEeqnarray*}{c}
\begin{array}{rl}
\tau_2&=\frac{T_{min}f_b-2\alpha T_2+\beta T_1+u_1}{f_u-f_{bw}}
\\
&=\frac{\frac{-\alpha T_2+\beta T_1+u_1}{T_2-T_{min}}-f_{bw}}{f_u-f_{bw}}(T_2-T_{min}).
\end{array}
\end{IEEEeqnarray*}
Since $\tau_2\in[0,T_2-T_{min}]$, if $T_{min}\leq T_2+\frac{\alpha T_2-\beta T_1-u_1}{f_{u}}$, the optimal solution exist and the solution is
\begin{equation*}
\begin{array}{rl}
F_2(t) = \left\{ \begin{array}{rl}
 f_{u} &\mbox{ if $t\in[t_2+T_{min},t_2+\tau_2]$} \\
 f_{bw} &\mbox{ if $t\in[t_2+\tau_2,t_2+T_2]$}
\end{array} \right.
\end{array}
\end{equation*}
Hence, by combining the results of $F_1$ and $F_2$, the optimal actuator force in $[t_2,t_2+T_2]$ which minimizes the total energy consumed for a single period is 
\begin{equation}
\begin{array}{rl}
F(t) = \left\{ \begin{array}{rl}
 f_{bw} &\mbox{ if $t\in[t_2, t_2+\tau_1]$} \\
 f_u &\mbox{ if $t\in[t_2+\tau_1,t_2+T_{min}+\tau_2]$}\\
 f_{bw} &\mbox{ if $t\in[t_2+\tau_2,t_2+T_2]$}.
\end{array} \right.
\end{array}
\label{eqn_optimal_G}
\end{equation}
In the next section, the necessary condition for existence of such a solution will be covered.

\subsubsection{STEP 4}
Here, we list the necessary conditions on $(u_1,v_1,T_{min})$ which gives the existence of a solution. Since $F(t)$ is a positive bounded function and $H,G$ and $I$ are increasing functions, the following bounds need to hold
\begin{IEEEeqnarray*}{c}
\int_{t_1}^{t_1+T_1}f_{bw}dt\leq H(t_1+T_1)\leq\int_{t_1}^{t_1+T_1}f_udt
\\
\int_{t_2}^{t_2+T_{min}}f_{bw}dt\leq G(t_2+T_{min})\leq\int_{t_2}^{t_2+T_{min}}f_{u}dt
\\
\int_{t_2+T_{min}}^{t_2+T_2}f_{bw}dt \leq G(t_2+T_2)-G(t_2+T_{min})
\\
\int_{t_2+T_{min}}^{t_2+T_2}f_{u}dt \geq G(t_2+T_2)-G(t_2+T_{min})
\\
\int_{t_3}^{t_3+T_3}f_{bw}dt\leq I(t_3+T_3)\leq\int_{t_3}^{t_3+T_3}f_{u}dt.
\end{IEEEeqnarray*}
Using the conditions for $T_{min}$ in STEP 3, the above bounds simplify to 
\begin{IEEEeqnarray}{c}
\beta(\frac{T(\eta+\rho)}{2(1+\rho)}-T_1)\geq u_1\geq\beta(\frac{T}{2}-T_1)
\label{eqn_bound_final_u}
\\
-K f_{bw} T_1-u_1\leq v_1\leq-2 f_{bw} T_1-u_1
\label{eqn_bound_final_v}
\\
\begin{array}{rl}
\max\{\gamma,T_2(1+\rho)-\gamma\eta\}\leq T_{min}
\\
T_{min}\leq \min\{\gamma\eta, T_2(1+\frac{\rho}{\eta})-\gamma\}
\end{array}
\label{eqn_bound_final_Tmin}
\\
0<T_1<T\frac{\rho}{1+\rho}
\label{eqn_bound_final_T1}
\\
1\leq \eta
\label{eqn_bound_final_eta}
\end{IEEEeqnarray}
where $\gamma=\frac{(\beta T_1+u_1)}{f_{u}}$, $K=\min\{1+\eta,(\eta-1)\frac{T_3}{T_1}-1\}$ and $\eta=\frac{f_u}{f_{bw}}\geq 1$. By choosing $T_1$ satisfying Eqn(\ref{eqn_bound_final_T1}), we found a set of admissible initial conditions, $(v_1,u_1)$, which satisfy Eqn(\ref{eqn_bound_final_u}) and Eqn(\ref{eqn_bound_final_v}). The bounds for $T_{min}$ are then well defined by the Eqn(\ref{eqn_bound_final_u}) and Eqn(\ref{eqn_bound_final_v}). 

\subsubsection{STEP 5} By substituting the optimal $G(t)$ into Eqn(\ref{eqn_work_final_L}), we can compute the power over unit distance from Eqn(\ref{eqn_P_U}). A long computation shows that Eqn(\ref{eqn_P_U}) ends up with a quadratic equation of $T_{min}$. By using the boundary condition for $T_{min}$, Eqn(\ref{eqn_bound_final_Tmin}), we can find the optimal $T_{min}$. Finally, by choosing $F$ and $I$ that satisfies Eqn(\ref{eqn_excursion_compute_final_L}), the first optimization problem is solved. 

\subsection{Solution to the second problem}
\subsubsection{Maximum average velocity} By using the boundary condition for $T_1$ in Eqn(\ref{eqn_bound_final_T1}) and the average velocity in Eqn(\ref{eqn_V_final}), we can compute $T_1$ which maximizes the average velocity. Since $\beta>0$ for all admissible $T_1$, we conclude that $V$ reaches its maximum when $T_1=T\frac{\rho}{1+\rho}$.

\subsubsection{Power over unit distance}
This problem is solved by STEP 5 of the first problem. Here, instead of $T_{min}$ that minimizes the performance, we find the optimal $(T_1,T_{min})$ pair for given boundary conditions. The boundary conditions for this problem are obtained from Eqn(\ref{eqn_bound_final_T1}) and Eqn(\ref{eqn_bound_final_Tmin}).

\section{Simulations}
\begin{figure}[!ht]
\centering
\includegraphics[scale=0.48]{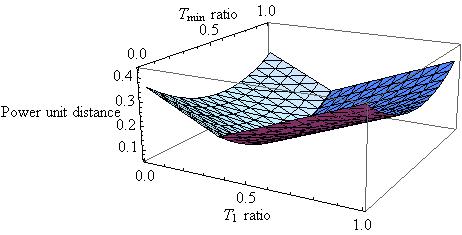}
\caption{A simulation for the power over an unit distance.}
\label{fig_sim_P_U}
\end{figure}

For simulation purpose, we choose the physical constraints as $T=10, f_{bw}=1, f_{fw}=0.1, f_u=5, d_1=40, t_1=0$ and $L=32.261$. Let $u_{ratio}$ and $v_{ratio}$ be the relative portion from the minimum to the maximum of its boundary, Eqn(\ref{eqn_bound_final_u}) and Eqn(\ref{eqn_bound_final_v}), respectively. Assume that $u_{ratio}=0.2$ and $v_{ratio}=0.5$. Now define $T_{1r}$ as the relative portion between minimum and maximum of the $T_1$ boundary. Similarly, define $T_{minr}$ as the relative portion between minimum and maximum of the $T_{min}$ boundary. Fig.\ref{fig_sim_P_U} shows the plot of the power over unit distance for all possible pairs $(T_{1r},T_{minr})$. 

\begin{figure}[!ht]
\centering
\includegraphics[width=2.9in]{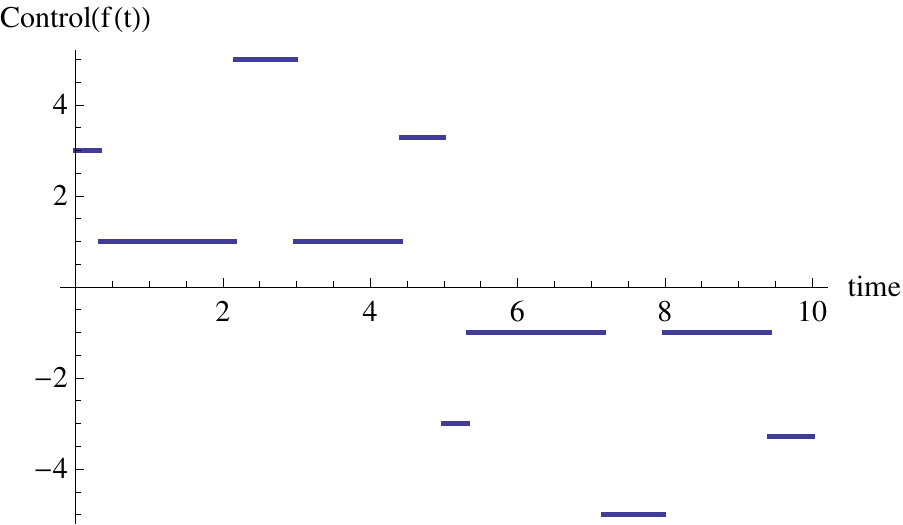}
\caption{The optimal actuator force model that minimize the power over unit distance.}
\label{fig_opt_f}
\end{figure}

The minimum occurs when $T_{1r}=0.363635$ and $T_{minr}=0.563214$. For a given $(T_1,T_{min})$, we can find the optimal actuator force model for $G$ by using Eqn(\ref{eqn_optimal_G}). Since we have freedom to choose $H$ and $I$ which satisfies the maximum excursion constraint, Eqn(\ref{eqn_excursion_compute_final_L}), one of the solutions that meet the maximum excursion constraint is shown in Fig.\ref{fig_opt_f}. 

\begin{figure}[!ht]
\centering
\includegraphics[width=2.9in]{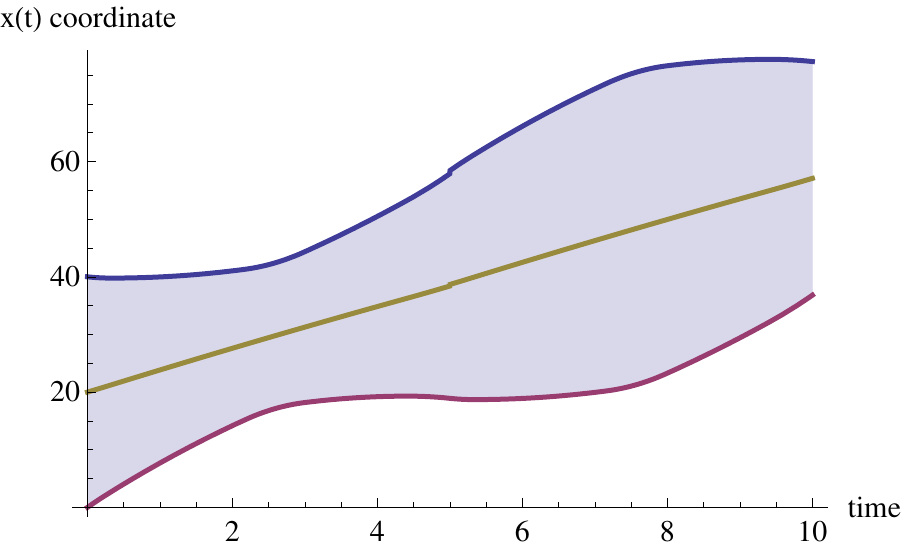}
\caption{A simulation for the motion in one period.}
\label{fig_sim_x}
\end{figure}

The minimum and maximum of excursion occur when $t=2.24$ and $t=7.74$, respectively.
Fig.\ref{fig_sim_x} shows the motion of each part. The blue(top) line is the trajectory of $x_2(t)$ and the red(bottom) line is the trajectory of $x_1(t)$. The middle line shows the motion of the center of mass. The optimal force uses its maximum in an interval about the point of minimum and maximum excursion.
\section{Conclusion}
Whereas previous work researched the admissible controls to be piecewise constant without a bound on the force, we have shown in this paper that such a shape of control can be also the locally optimal for the case when we optimize over the piecewise continuous and bounded controls. However, due to the freedom in choosing $H$ and $I$, other solutions having the same performance exist. The result shows that it is necessary to allocate the maximum allowable force before and after the size of the worm gets close to its minimum or maximum length. 

\section{Discussion}
Since the set, $\{t\in[t_2,t_2+T_2]:\lambda_{i}(t)=0\}$, only contains two points, there does not exist a singular solution which performs better than the local optimal solution. In addition, all the admissible controls $F$ in Case1 and Case3, which satisfy the equality constraint Eqn(\ref{eqn_excursion_compute_final_L}), generate the same power performance. If we consider the problem when there exist an additional penalty before and after the swtich happens, then one may find the unique optimal control among the solutions to Eqn(\ref{eqn_excursion_compute_final_L}).



\bibliographystyle{IEEEtran}
\bibliography{mtns2014}

\end{document}